%Differs from 3-10-98 version by implementing referee's changes

%\magnification=1200

%%%%%%%%%%%%%%%%%%%%%%%%%%%%%%% Macro Definitions %%%%%%%%%%%%%%%%%%%%%%%%%%%

\def\title#1{{\titlefont\noindent #1\bigskip}}

\def\author#1{{\largefont\noindent #1}\medskip}

\def\beginlinemode{\endmode
 \begingroup\obeylines\def\endmode{\par\endgroup}}
\let\endmode=\par

\newbox\theaddress
\def\address{\smallskip\beginlinemode\parindent 0in\getaddress}
{\obeylines
\gdef\getaddress #1 
 #2
 {#1\gdef\addressee{#2}%
   \global\setbox\theaddress=\vbox\bgroup\raggedright%
    \everypar{\hangindent2em}#2
   \def\endaddress{\egroup\endgroup \copy\theaddress \medskip}}}

\def\thanks#1{\footnote{}{\eightpoint #1}}

\long\def\Abstract#1{{\eightpoint\narrower\vskip\baselineskip\noindent
#1\smallskip}}

\def\skipfirstword#1 {}

\def\ir#1{\csname #1\endcsname}

\newdimen\currentht
\newbox\droppedletter
\newdimen\droppedletterwdth
\newdimen\drophtinpts
\newdimen\dropindent

\def\irrnSection#1#2{\edef\tttempcs{\ir{#2}}
\currentht-\pagetotal\advance\currentht by-\ht\footins
\advance\currentht by\vsize
\ifdim\currentht<1.5in\par\vfill\eject\else\vbox to.25in{\vfil}\fi
{\largefont\noindent{}\expandafter\skipfirstword\tttempcs. #1}
\vskip\baselineskip }

\def\irSubsection#1#2{\edef\tttempcs{\ir{#2}}
\vskip\baselineskip\penalty-3000
{\bf\noindent{}\expandafter\skipfirstword\tttempcs. #1}
\vskip\baselineskip}

\def\irSubsubsection#1#2{\edef\tttempcs{\ir{#2}}
\vskip\baselineskip\penalty-3000
{\bf\noindent \expandafter\skipfirstword\tttempcs. #1}
\vskip\baselineskip}

\def\References{\vbox to.25in{\vfil}\noindent{}{\bf References}
\vskip\baselineskip\par}

\def\baselinebreak{\par \ifdim\lastskip<\baselineskip
         \removelastskip\penalty-200\vskip\baselineskip\fi}

\long\def\prclm#1#2#3{\baselinebreak
\noindent{\bf \csname #2\endcsname}:\enspace{\sl #3\par}\baselinebreak}

\def\Prf{\noindent{\bf Proof}: }

\def\rem#1#2{\baselinebreak\noindent{\bf \csname #2\endcsname}:\enspace }

\def\qed{{\hfill$\diamondsuit$}\vskip\baselineskip}

\def\bibitem#1{\par\indent\llap{\rlap{\bf [#1]}\indent}\indent\hangindent
2\parindent\ignorespaces}

\long\def\eatit#1{}

\def\leftheadlinetext{}
\def\rightheadlinetext{}

\def\leftheadline{{\eightrm\folio\hfil \leftheadlinetext\hfil}}
\def\rightheadline{{\eightrm\hfil\rightheadlinetext\hfil\folio}}

\headline={\ifnum\pageno=1\hfil\else
\ifodd\pageno\rightheadline\else\leftheadline\fi\fi}

\def\tenpoint{\def\rm{\fam0\tenrm}
\textfont0=\tenrm \scriptfont0=\sevenrm \scriptscriptfont0=\fiverm
\textfont1=\teni \scriptfont1=\seveni \scriptscriptfont1=\fivei
\def\mit{\fam1} \def\oldstyle{\fam1\teni}
\textfont2=\tensy \scriptfont2=\sevensy \scriptscriptfont2=\fivesy
\def\cal{\fam2}
\textfont3=\tenex \scriptfont3=\tenex \scriptscriptfont3=\tenex
\def\it{\fam\itfam\tenit} % \it is family 4
\textfont\itfam=\tenit
\def\sl{\fam\slfam\tensl} % \sl is family 5
\textfont\slfam=\tensl
\def\bf{\fam\bffam\tenbf} % \bf is family 6
\textfont\bffam=\tenbf \scriptfont\bffam=\sevenbf
\scriptscriptfont\bffam=\fivebf
\def\tt{\fam\ttfam\tentt} % \tt is family 7
\textfont\ttfam=\tentt
\normalbaselineskip=12pt
\setbox\strutbox=\hbox{\vrule height8.5pt depth3.5pt  width0pt}%
\normalbaselines\rm}

\def\eightpoint{\def\rm{\fam0\eightrm}%
\textfont0=\eightrm \scriptfont0=\sixrm \scriptscriptfont0=\fiverm
\textfont1=\eighti \scriptfont1=\sixi \scriptscriptfont1=\fivei
\def\mit{\fam1} \def\oldstyle{\fam1\eighti}%
\textfont2=\eightsy \scriptfont2=\sixsy \scriptscriptfont2=\fivesy
\def\cal{\fam2}%
\textfont3=\tenex \scriptfont3=\tenex \scriptscriptfont3=\tenex
\def\it{\fam\itfam\eightit} % \it is family 4
\textfont\itfam=\eightit
\def\sl{\fam\slfam\eightsl} % \sl is family 5
\textfont\slfam=\eightsl
\def\bf{\fam\bffam\eightbf} % \bf is family 6
\textfont\bffam=\eightbf \scriptfont\bffam=\sixbf
\scriptscriptfont\bffam=\fivebf
\def\tt{\fam\ttfam\eighttt} % \tt is family 7
\textfont\ttfam=\eighttt
\normalbaselineskip=9pt%
\setbox\strutbox=\hbox{\vrule height7pt depth2pt  width0pt}%
\normalbaselines\rm}

\def\largefont{\def\rm{\fam0\largerm}
\textfont0=\largerm \scriptfont0=\largescriptrm \scriptscriptfont0=\tenrm
\textfont1=\largei \scriptfont1=\largescripti \scriptscriptfont1=\teni
\def\mit{\fam1} \def\oldstyle{\fam1\teni}
\textfont2=\largesy %\scriptfont2=\sevensy \scriptscriptfont2=\fivesy
\def\cal{\fam2}
%\textfont3=\largeex %\scriptfont3=\tenex \scriptscriptfont3=\tenex
\def\it{\fam\itfam\largeit} % \it is family 4
\textfont\itfam=\largeit
\def\sl{\fam\slfam\largesl} % \sl is family 5
\textfont\slfam=\largesl
\def\bf{\fam\bffam\largebf} % \bf is family 6
\textfont\bffam=\largebf %\scriptfont\bffam=\sevenbf \scriptscriptfont\bffam=\fivebf
\def\tt{\fam\ttfam\largett} % \tt is family 7
\textfont\ttfam=\largett
\normalbaselineskip=17.28pt
\setbox\strutbox=\hbox{\vrule height12.25pt depth5pt  width0pt}%
\normalbaselines\rm}

\def\titlefont{\def\rm{\fam0\titlerm}
\textfont0=\titlerm \scriptfont0=\largescriptrm \scriptscriptfont0=\tenrm
\textfont1=\titlei \scriptfont1=\largescripti \scriptscriptfont1=\teni
\def\mit{\fam1} \def\oldstyle{\fam1\teni}
\textfont2=\titlesy %\scriptfont2=\sevensy \scriptscriptfont2=\fivesy
\def\cal{\fam2}
%\textfont3=\largeex %\scriptfont3=\tenex \scriptscriptfont3=\tenex
\def\it{\fam\itfam\titleit} % \it is family 4
\textfont\itfam=\titleit
\def\sl{\fam\slfam\titlesl} % \sl is family 5
\textfont\slfam=\titlesl
\def\bf{\fam\bffam\titlebf} % \bf is family 6
\textfont\bffam=\titlebf %\scriptfont\bffam=\sevenbf \scriptscriptfont\bffam=\fivebf
\def\tt{\fam\ttfam\titlett} % \tt is family 7
\textfont\ttfam=\titlett
\normalbaselineskip=24.8832pt
\setbox\strutbox=\hbox{\vrule height12.25pt depth5pt  width0pt}%
\normalbaselines\rm}

\nopagenumbers

\font\eightrm=cmr8
\font\eighti=cmmi8
\font\eightsy=cmsy8
\font\eightbf=cmbx8
\font\eighttt=cmtt8
\font\eightit=cmti8
\font\eightsl=cmsl8
\font\sixrm=cmr6
\font\sixi=cmmi6
\font\sixsy=cmsy6
\font\sixbf=cmbx6

\font\largerm=cmr12 at 17.28pt
\font\largei=cmmi12 at 17.28pt
\font\largescriptrm=cmr12 at 14.4pt
\font\largescripti=cmmi12 at 14.4pt
\font\largesy=cmsy10 at 17.28pt
\font\largebf=cmbx12 at 17.28pt
\font\largett=cmtt12 at 17.28pt
\font\largeit=cmti12 at 17.28pt
\font\largesl=cmsl12 at 17.28pt

\font\titlerm=cmr12 at 24.8832pt
\font\titlei=cmmi12 at 24.8832pt
\font\titlesy=cmsy10 at 24.8832pt
\font\titlebf=cmbx12 at 24.8832pt
\font\titlett=cmtt12 at 24.8832pt
\font\titleit=cmti12 at 24.8832pt
\font\titlesl=cmsl12 at 24.8832pt

\tenpoint

%%%%%%%%%%%%%%%%%%%%%%%%%%%%%%% Internal References %%%%%%%%%%%%%%%%%%%%%%%%%%%
%Internal references are implemented via the \ir command,
%hence \ir{intro} rather than \intro. This makes it easier
%to search for all internal references in the text.

%%%%%%%%%%%%%%%%%%%%%%%%%%  BEGIN PAPER  %%%%%%%%%%%%%%%%%%%%%%%%%%%%%%%%
\def\manyby{\hbox to.75in{\hrulefill}}
\hsize 6.5in 
\vsize 9.2in

\tolerance 3000
\hbadness 3000

\def\item#1{\par\indent\indent\llap{\rlap{#1}\indent}\hangindent
2\parindent\ignorespaces}

\def\itemitem#1{\par\indent\indent
\indent\llap{\rlap{#1}\indent}\hangindent
3\parindent\ignorespaces}

\def\itemitemitem#1{\par\indent\indent\indent
\indent\llap{\rlap{#1}\indent}\hangindent
4\parindent\ignorespaces}

\def\itemitemitemitem#1{\par\indent\indent\indent\indent
\indent\llap{\rlap{#1}\indent}\hangindent
5\parindent\ignorespaces}

\def\refCat{1}
\def\refFi{2}
\def\refGO{3}
\def\refGGR{4}
\def\trans{5}

\def\vanc{6}
\def\ravello{7}
\def\mtnwest{8}
\def\antican{9}
\def\fatpts{10}
\def\newfatpts{11}

\def\refHa{12}
\def\refHi{13}
\def\refMa{14}
\def\refMu{15}
\def\refNtwo{16}

\def\C#1{\hbox{$\cal #1$}}

\def\div#1 #2 #3 #4 #5 #6 #7 #8 {#1L #2E_1 #3E_2 #4E_3 #5E_4 #6E_5 #7E_6 #8E_7}

\def\pr#1{\hbox{{\bf P}${}^{#1}$}}
\def\leftheadlinetext{Brian Harbourne}
\def\rightheadlinetext{Fat Point Algorithm}

\title{An Algorithm for Fat Points on \pr2}

\author{Brian Harbourne}

\address
Department of Mathematics and Statistics
University of Nebraska-Lincoln
Lincoln, NE 68588-0323
email: bharbourne@unl.edu
WEB: http://www.math.unl.edu/$\sim$bharbour/
\smallskip
September 13, 1999 \endaddress
\vskip-\baselineskip

\thanks{\vskip -6pt
\noindent This work benefitted from a National Science Foundation grant.
I also thank the referee for his suggestions and the resulting
improvement in the exposition.
\smallskip
\noindent 1991 {\it Mathematics Subject Classification. } 
Primary 13P10, 14C99. 
Secondary 13D02, 13H15.
\smallskip
\noindent {\it Key words and phrases. }  Generators, syzygies,
resolution, fat points, maximal rank, plane, Weyl group.\smallskip}

\vskip\baselineskip
\Abstract{Abstract: Let $F$ be a
divisor on the blow-up $X$ of \pr2 at $r$ general points
$p_1, \ldots, p_r$ and let $L$ be the total transform of 
a line on \pr2. An approach is presented for reducing 
the computation of the
dimension of the cokernel of the natural map
$\mu_F:\Gamma(\C O_X(F))\otimes\Gamma(\C O_X(L))\to
\Gamma(\C O_X(F)\otimes\C O_X(L))$ to the case that
$F$ is ample. As an application, a formula
for the dimension of the cokernel of $\mu_F$ is
obtained when $r= 7$, completely solving the problem of
determining the modules in minimal free 
resolutions of fat point subschemes
$m_1p_1+\cdots+ m_7p_7\subset \pr2$.
All results hold for an arbitrary algebraically closed 
ground field $k$.}
\vskip\baselineskip

\irrnSection{Introduction}{intro}
Let $p_1, \ldots, p_r\in \pr N$ be general points in projective space,
let $m_1,\ldots,m_r$ be nonnegative integers and let
$I(p_i)$ be the homogeneous ideal (in the homogeneous 
coordinate ring $R= k[\pr N]$ of \pr N) generated by all homogeneous
polynomials vanishing at $p_i$.
A {\it fat point\/} subscheme
$Z=m_1p_1+\cdots+ m_rp_r\subset
\pr N$ is the subscheme  
corresponding to the homogeneous ideal $I(Z)=
I(p_1)^{m_1}\cap\cdots\cap I(p_r)^{m_r}$ (which it is easy to see is
generated by all homogeneous 
polynomials vanishing at each point $p_i$ to order
at least $m_i$). If $m_i\le1$ for all $i$, we say 
$Z$ is a {\it thin point} subscheme.
We denote by $I(Z)_t$ the homogeneous component of $I(Z)$ of degree $t$.

The first module in any minimal free homogeneous resolution of $I(Z)$
is, up to graded isomorphism, $\oplus_t R[-t]^{\nu_t}$, where
$\nu_t$ (or $\nu_t(Z)$ if for clarity $Z$ needs to be specified)
is the dimension of the cokernel of the obvious
multiplication map $\mu_{t-1}(Z): I(Z)_{t-1}\otimes R_1\to I(Z)_t$.
(More concretely, $\nu_t$ is the number of generators in degree $t$ of
any minimal set of homogeneous generators of $I(Z)$.)

In the case of a thin point subscheme
$Z=p_1+\cdots+ p_r$ (with $p_i$ general), 
the dimensions of the homogeneous components $I(Z)_j$ 
are known so one can determine $\nu_t(Z)$ from the rank of $\mu_{t-1}(Z)$,
and the maximal rank conjecture of [\refGO, \refGGR] is that
$\mu_t$ should be of maximal rank for all $t$ (meaning that 
$\mu_t$ should always be either injective or surjective).
Although this conjecture has been verified in a number of cases
(including $N=2$), it remains open in general.
For the more general but analogous situation of fat points, 
no conjecture has been put forward. 
This is partly because
the multiplication maps often fail to have maximal rank, 
and partly because little is known about how otherwise
the ranks and numbers of generators should behave, but 
also because one typically first wants to understand
Hilbert functions, and Hilbert functions 
of fat point ideals are themselves not yet well understood.

However, understanding of Hilbert functions for $N=2$, 
although not complete, is much better
than in higher dimensions. Indeed, there are comprehensive conjectures 
(see [\vanc, \ravello, \refHi]) which 
in various situations are known to hold. Thus some attention has
begun to be paid to the behavior of
generators and resolutions of ideals of fat point subschemes for $N=2$, 
both for its own interest and as an initial means 
of developing one's understanding in general.

So, for the rest of this paper we will assume $N=2$,
in which case, since a fat points subscheme of \pr2 is
arithmetically Cohen-Macaulay, 
a minimal free graded resolution of $I(Z)$ is of the form
$0\to F_1\to F_0\to I(Z)\to 0$; the values $\nu_t$ determine
$F_0=\oplus_t R[-t]^{\nu_t}$, which with the Hilbert function of $I(Z)$
then determines $F_1$. Thus, for $N=2$, given the numbers of generators
and the Hilbert function of $I(Z)$, one also has
the modules in a minimal free resolution of $I(Z)$.

\irSubsection{The Particular Interest of $\nu_{\beta+1}$ and $r=7$}{degbeta}
Denote by $\alpha(Z)$ (or by just
$\alpha$ when $Z$ is understood) the 
least degree $t$ such that $I(Z)_t\ne 0$ 
and by $\beta(Z)$ the least degree $t$ such that the base locus
of $I(Z)_t$ is 0-dimensional (said alternately, 
$\beta(Z)$ is the least degree $t$ such that the elements of $I(Z)_t$
have no nontrivial common divisor). Given $Z=m_1p_1+\cdots+m_rp_r$,
if the points $p_i$ are sufficiently general, conjecturally
(see [\vanc, \ravello, \refHi]) the regularity of 
$I(Z)$ is at most $\beta(Z)+1$, assuming which
the general problem of finding $\nu_t$ 
reduces by Lemmas 2.9 and 2.10 of [\fatpts] 
to computing Hilbert functions and $\nu_{\beta+1}$, and thus 
the case $t=\beta+1$ is of particular interest.
In [\refFi] Fitchett develops
a means of handling $\nu_{\beta+1}$ 
in the case that $\alpha<\beta$, but
it remains unclear what to do when $\alpha=\beta$.
The naive hope that in this situation
$\mu_\beta$ might have maximal rank
is quashed by examples from [\newfatpts] showing
that maximal rank can fail.

It is for $r=7$ general points of \pr2
that we can first hope to begin to understand 
the source of such failures, since for $r\le 6$ there are none.
For example, for $r\le 5$, it follows from [\refCat] 
that $\nu_{\beta+1}=0$ always holds (or see [\fatpts]), 
and hence that $\mu_\beta$
has maximal rank, and for $r=6$, 
although $\nu_{\beta+1}$ need not always vanish, [\refFi] shows 
that $\mu_\beta$ always has maximal rank.

\irSubsection{The Geometric Translation}{geomtrans}
Thus in this paper we resolve the problem
of determining $\nu_{\beta+1}$ when $r=7$,
thereby working out the resolution of ideals
defining fat point subschemes 
involving $r=7$ general points of \pr2, taking,
as did [\refFi, \fatpts, \newfatpts],
a geometric approach in which we obtain
results for line bundles on certain rational surfaces. 
Those readers unfamiliar with
the by-now standard translation of questions about fat
points in \pr2 to questions about line bundles on
blow ups of \pr2 may find it helpful to refer to
[\fatpts]. In particular, for any fat point subscheme
$Z=m_1p_1+\cdots+m_rp_r\subset\pr2$,
there is for each degree $t$ a corresponding divisor $F$
(which is effective and 
numerically effective when $t\ge\beta$)
on the blow up $X$ of \pr2 at the points $p_i$ such that 
the dimension of $I(Z)_t$ is $h^0(X,\C O_X(F))$,
and $\nu_{t+1}$ is the dimension  
of the cokernel $\C S(F,L)$ of the natural map
$\mu_F:\Gamma(\C O_X(F))\otimes\Gamma(\C O_X(L))\to
\Gamma(\C O_X(F)\otimes\C O_X(L))$, where $L\subset X$ 
is the total transform to $X$ of a line on \pr2. 

For the reader's convenience, we recall some standard notions
from geometry. A divisor (as always, on a given 
smooth projective surface) 
which is a nonnegative (integer) linear
combination of curves is said to be {\it effective}.
A divisor $F$ (or its linear equivalence class $[F]$) 
is {\it numerically effective\/} if 
$F\cdot C\ge 0$ for every effective divisor $C$, while
an {\it ample\/} divisor is one
whose intersection with every effective divisor 
is positive. On the other hand, an {\it exceptional curve\/}
is a smooth rational curve of 
self-intersection $-1$; for example, 
the curve obtained by blowing up a smooth point on a 
projective surface is an exceptional curve. The exceptional curves
on a smooth projective rational surface are known;
see [\refNtwo, \refMa].

Now let $F$ be a 
divisor on a surface $X$ obtained by blowing up distinct points $p_1,
\ldots, p_r$ of \pr2. 
Let $L$ be the total transform to $X$ of a line on \pr2. 
We denote by $\hbox{Cl}(X)$
the group of divisors on $X$ modulo linear equivalence. This quotient, 
the {\it divisor class group}, is a free abelian group. 
The classes $[L]$, $[E_1]$, $\ldots$, $[E_r]$ (where $E_i$ is 
the exceptional curve obtained by blowing up $p_i$) give a basis of
$\hbox{Cl}(X)$ which we refer to as an {\it exceptional configuration}.
(In this notation, the divisor $F$ referred to above, 
corresponding in degree $t$ to $Z=m_1p_1+\cdots+m_rp_r\subset\pr2$,
is $F=tL-m_1E_1-\cdots-m_rE_r$.)
Also, there is a bilinear form,
the intersection form, on $\hbox{Cl}(X)$, in which the
basis elements $[L]$, $[E_1]$, $\ldots$, 
$[E_r]$ are orthogonal and such that
$-[L]^2=[E_1]^2= \ldots=[E_r]^2=-1$.

\irSubsection{Discussion of Results}{resdis}
So, in fact, in this paper we solve the problem of
computing the dimension of the cokernel of
$\mu_F$ for arbitrary divisors $F$ on
a blow up $X$ of \pr2 at 7 general points.
Our solution is, first, algorithmically to reduce
to the case that $F$ is ample, and second,
to show that $\mu_F$ is surjective when $F$ is ample.
This approach mimics what was already known
concerning the determination of the Hilbert function of
$I(Z)$ involving $r\le 9$ general points of \pr2. 
As mentioned above, given $t$, there is a corresponding 
divisor $F$ on the blow up $X$ of \pr2 at the $r$ points
such that the dimension
of $I(Z)_t$ is equal to $h^0(X,\C O_X(F))$.
But as shown in [\trans] and [\antican] 
and as is discussed in [\mtnwest],
one can for $r\le 9$ general points algorithmically reduce
the computation of $h^0(X,\C O_X(F))$
to the case that $F$ is numerically effective,
in which case $h^0(X,\C O_X(F))=(F^2-F\cdot K_X)/2+1$.

For computing the dimension of the cokernel of
$\mu_F$, the reduction to ample $F$ depends on 
three hypotheses, (A1), (A2) and (A3),
which we explicitly mention below and 
which are known to hold for a divisor on a 
surface obtained by blowing up
$r\le 8$ general points. If we consider 
$r$ generic points, these hypotheses
continue to hold for $r=9$ and they
have been conjectured to hold for all $r$.
(Alternatively, since (A2) and (A3) 
pose conditions on possibly infinite sets of divisors,
only finitely many of which are relevant
at any one time, (A2) and (A3) are slightly stronger than 
needed for our purposes. So in fact
slightly weaker but more complicated versions
of (A2) and (A3) can also be used
which are known to hold for $r=9$ general points,
and which are conjectured to hold for any $r$ general points.)

Thus the main difficulty of working out resolutions
of fat point subschemes involving $r>7$ general
points is not the reduction to ampleness.
It is rather that for $r>7$, 
$\mu_F$ need not be surjective when $F$ is ample.
If $r<7$, the surjectivity of $\mu_F$ for an ample divisor
$F$ already follows from [\refCat, \refFi, \fatpts]. In this paper,
\ir{amplesthm} extends this to $r=7$. 
But in both cases, the tools used to show
that the cokernel of $\mu_F$ vanishes
for an ample divisor $F$ only bound
the dimension of the cokernel in general, 
and the bounds obtained are not always delicate enough to
pin down the dimension of the cokernel completely
when $r>7$. If the set of problematical cases
were not too large, one could hope to handle
them ad hoc, and this seems possible
in case $r=8$, but as of this writing this 
does not seem workable for $r>8$. (This is related
to the fact that $K_X^\perp$ is negative definite for $r<9$,
but indefinite for all $r>8$, where 
$K_X$ denotes the canonical class of $X$ and
$K_X^\perp$ denotes the 
subspace of $\hbox{Cl}(X)$
of classes of divisors $D$ with $D\cdot K_X=0$.)

Our main result is \ir{mainthm}
(but see also \ir{lastcor}),
which explicitly determines
the dimension of the cokernel
of $\mu_F$ for any numerically effective divisor
$F$ on $X$. We regard this as our main result because,
from the point of view of a homogeneous 
ideal defining a fat points subscheme 
$Z=\sum_ip_i\subset\pr2$, numerically 
effective divisors are more natural than
ample divisors. For example, if $I(Z)_t\ne 0$,
let $V\subset R$ be the subspace
of elements obtained from $I(Z)_t$
by dividing out by a greatest divisor
common to all of the elements of $I(Z)_t$.
Then, under the standard translation, $V$ corresponds
to $|F|$ for some numerically effective divisor $F$
on the blow up of \pr2 at the points $p_i$.
From the dimensions of the linear system $|F|$
and of the cokernel of $\mu_F$ we can find the dimension
of the kernel of $\mu_F$, which is the same as that of
$V\otimes R_1\to R$, whose kernel has the same dimension 
as that of $\mu_t(Z)$, which with
the dimensions of $I(Z)_t$ and $I(Z)_{t+1}$,
allows us to compute the dimension of the cokernel
of $\mu_t(Z)$, and hence $\nu_{t+1}(Z)$.

Putting it all together, we
obtain an algorithm for determining 
$\nu_t(Z)$ for each $t$ for any
fat point subscheme $Z$ involving $r\le 7$
general points of \pr2. As discussed above, since the Hilbert
function of $I(Z)$ is known, this
gives an algorithm for determining up to graded isomorphism
the modules in the minimal free 
resolutions of the ideal $I(Z)$. (As of this writing,
an implementation of this algorithm can be run via the
World-Wide Web at the author's web site;
the specific web address is 
\hbox{http://www.math.unl.edu/$\sim$bharbour/cgi-bin/7fatpts.cgi .})

\irSubsection{Algorithm's Underlying Assumptions}{algass}
Our algorithm assumes that:
\item{$\bullet$}(A1) $X$ is obtained by blowing 
up $r$ distinct points of \pr2, that 
\item{$\bullet$}(A2) the only curves of negative 
self-intersection on $X$ are exceptional curves, 
and that 
\item{$\bullet$}(A3) $h^1(X, \C O_X(F))=0$ for any effective,
numerically effective divisor $F$.

By [\mtnwest], (A3) holds for any $r\le 8$ points,
general or not. For $r\le 8$ general points, $X$ is Del Pezzo, so
$-K_X$ is ample, so by adjunction
(A2) holds too. 

More generally, since $K_X=-3[L]+[E_1+\cdots+E_r]$,
it follows for all $r$ by duality that 
$h^2(X, \C O_X(F))=0$ whenever $F\cdot L>-3$ (such as is 
the case if $F$ is numerically effective or effective). In addition,
[\mtnwest] shows $F^2\ge 0$ for any numerically effective divisor $F$.
Moreover, for an arbitrary divisor $F$ it is true that
$h^0(X, \C O_X(F))=h^0(X, \C O_X(F-E))$
if $E$ is effective, reduced and irreducible with $F\cdot E<0$.
By iteratively replacing $F$ by $F-E$ whenever $E$ is
an exceptional curve and $F\cdot E<0$ (see \ir{Wrem}), we 
thus eventually obtain a divisor
$F$ such that either $F\cdot L<0$, and hence
$h^0(X, \C O_X(F))=0$, or such that $F\cdot L\ge 0$ and
$F\cdot E\ge 0$ for every exceptional curve $E$.

But in the latter case, $h^0(X, \C O_X(F))>0$
if and only if $(F^2-K_X\cdot F)/2+1>0$.
This is because if $h^0(X, \C O_X(F))>0$, then
$F$ would be numerically effective (because 
$F$ meets all exceptional curves nonnegatively
and by (A2) there are no other curves of 
negative self-intersection), and hence by (A3)
$h^0(X, \C O_X(F))=(F^2-K_X\cdot F)/2+1$.
Conversely, $F\cdot L\ge 0$ means $h^2(X, \C O_X(F))=0$
so $(F^2-K_X\cdot F)/2+1> 0$ implies $h^0(X, \C O_X(F))>0$ 
by Riemann-Roch (and then, as before,
$h^0(X, \C O_X(F))=(F^2-K_X\cdot F)/2+1$).

In any case, we end up knowing $h^0(X, \C O_X(F))$.
Thus it is not an additional assumption
to assume that $h^0(X, \C O_X(F))$ is 
always available for any divisor class $[F]$ on $X$. 

But since we can assume that we always can compute
$h^0(X, \C O_X(F))$, we can also assume 
that given the class of any effective
divisor $F$, we can determine the classes of the divisors occurring
as fixed components of $|F|$. This is because
for any effective divisor $F$ there is
in terms of the classes $[L]$, $[E_1]$, $\ldots$, $[E_r]$ 
a finite list
of classes such that if $C$ is an effective divisor for which
$F-C$ is effective, then the class of $C$ is on the list:
if $[F]=d[L]-\sum_{i>0}m_i[E_i]$ and $[C]=d'[L]-\sum_{i>0}m_i'[E_i]$,
where $F$, $C$ and $F-C$ are all effective,
then $d\ge d'\ge0$, $d'\ge m_i'$ and $d-d'-m_i\ge -m_i'$ for all $i>0$.
Thus we have only finitely many
classes to test, the test being that 
$C$ is a fixed component of $|F|$ if and only if
$h^0(X, \C O_X(F))=h^0(X, \C O_X(F-C))$. 

We also have:

\prclm{Lemma}{introlem}{Given (A1), (A2) and (A3), 
let $F\ne 0$ be an effective divisor on $X$
with $F\cdot E>0$ for every exceptional curve
$E$ and such that $|F|$ is fixed component free.
Then either $F^2>0$ and $F$ is ample, or
$F^2=0$ and $|F|$ is composed with a pencil $|D|$
where $D$ is a smooth rational curve.}

\Prf Note that $F$ is numerically effective.
If $F^2>0$, by the Hodge Index Theorem and (A2),
there can be no effective divisors $C$ with $F\cdot C=0$,
so $F$ is ample. 

If, instead, $F^2=0$, we want to show that
$|F|=|tD|$ for some smooth rational curve $D$. 
By Riemann-Roch and (A3),
$1<h^0(X,\C O_X(F))=1-K_X\cdot F/2$, so $-K_X\cdot F>0$.
In any case, $|F|$ defines a morphism 
to \pr1. By Stein factorization,
$F$ is linearly equivalent to $tG$ for some $t>0$, where $|G|$ defines
a morphism to \pr1 with connected fibers, $G^2=0$ and $G$ is free. 
If $G$ is irreducible, then $G=sC$ for some prime divisor $C$,
so $C^2=0$ and $-K_X\cdot C>0$ (since $-K_X\cdot F>0$), hence by adjunction 
$C$ is a smooth rational curve and $|G|$ and hence $|F|$ is composed 
with the pencil $|C|$. If $G$ is not irreducible, 
then among the components of $G$ are two distinct
reduced and irreducible components $B$ and $C$ which meet.
I.e., $B\cdot C>0$, and since $G^2=0$, we also 
have $B^2<0$ and $C^2<0$,
so, by (A2), $B$ and $C$ are exceptional
and, since $(B+C)\cdot G=0$, also $B\cdot C=1$. 
By (A3) it follows that $|B+C|$
is a pencil, all elements of which are either
irreducible or sums of exceptional curves. As the 
latter can happen in only finitely many ways, a general
element $D$ of $|B+C|$ is irreducible, hence as before,
a smooth rational curve. Moreover, $G-B-C$ is effective
and $D+(G-B-C)$ is an element of $|G|$ and thus connected,
but $D+(G-B-C)\ne G$ so $D$ is disjoint from
$G$; hence $G-B-C=0$ so $G$ is linearly equivalent to $D$
and $|F|$ is composed with the pencil $|D|$.\qed

\rem{Remark}{Wrem} The procedure described
above for determining $h^0(X,\C O_X(F))$
depends on checking $F\cdot E$ for all
exceptional curves $E$. This is no problem
when $r\le 8$, since then there are only 
finitely many exceptional curves. More generally,
there is an action on $\hbox{Cl}(X)$ by a Weyl group, $W$;
see \ir{sbkgrnd}. Given (A1) and (A2),
in the terminology of
the proof of Theorem 2.1 of [\trans], no nodal classes are
effective, so $W$ acts transitively on the exceptional configurations
of $X$, and, for $r>2$, 
the set of classes of exceptional curves is precisely
a single $W$-orbit. For simplicity, let us say $r>2$.
Then the proof of Theorem 2.1 [\trans]
gives an algorithm (under different hypotheses
but still applicable here) for finding an element $w\in W$ such that
either: (i) $w(F)\cdot L<0$; (ii) $w(F)\cdot E_i<0$ for some $i>0$;
or (iii) $w(F)$ is a nonnegative sum of the classes
$[L], [L-E_1], [2L-E_1-E_2], [3L-E_1-E_2-E_3],\ldots,[3L-E_1-\cdots-E_r]$.
But $w(F)\cdot L<0$ means that $w(F)\cdot [C]<0$,
where $C$ is either $E_i$ for some $i>0$
or $L-E_1-E_2$, so in cases (i) or (ii) with $C$ being 
one of the exceptional curves whose class is
$[E_i]$ for some $i>0$ or $[L-E_1-E_2]$,
we have $F\cdot w^{-1}C<0$. In case (iii),
it is easy to check that the classes $[L], [L-E_1]$ and $[2L-E_1-E_2]$
are numerically effective while $[3L-E_1-\cdots-E_i]$ meets every 
exceptional class nonnegatively (since $[3L-E_1-\cdots-E_i]
=-K_X+[E_{i+1}+\cdots+E_r]$), and hence that $F\cdot E\ge 0$ 
for every exceptional curve $E$. Thus, for an arbitrary $F$,
we have an effective means of finding an $E$ with $F\cdot E<0$, or of deciding
none such exists.
\qed

\irSubsection{The Algorithm}{thealg}
So here is our algorithm. Assume (A1), (A2) and (A3)
and let $[F]$ be a 
divisor class on $X$. Our goal is to compute
$\hbox{dim}\,\hbox{ker}(\mu_F)$, from which we can obtain
our ultimate goal of computing $\hbox{dim}\,\hbox{cok}(\mu_F)$
via the obvious formula $\hbox{dim}\,\hbox{cok}(\mu_F)=h^0(X,\C O_X(F+L))
-3h^0(X,\C O_X(F))+\hbox{dim}\,\hbox{ker}(\mu_F)$.
The following algorithm reduces the problem of
computing $\hbox{dim}\,\hbox{ker}(\mu_F)$ for an arbitrary $F$ to
the case that $F$ is ample.

\vskip\baselineskip
\noindent {\bf START}: given $F$, compute $h^0(X, \C O_X(F))$.
\item{I.} If $h^0(X, \C O_X(F))\le 1$, then
clearly $\mu_F$ is injective: {\bf STOP}.
\item{II.} Assume $h^0(X, \C O_X(F))>1$:
\itemitem{1.} If $|F|$ has a fixed component $C$,
then clearly $\mu_F$ and $\mu_{F-C}$ have kernels of the same dimension,
and we replace $F$ by $F-C$. After a finite number 
of such subtractions, we reduce to the case
that $F$ is effective and $|F|$ is fixed component free,
without changing $\hbox{dim}\,\hbox{ker}(\mu_F)$: go to step 2.
\itemitem{2.} Assume $h^0(X, \C O_X(F))\ge 2$ and $|F|$ has no fixed components.
\itemitemitem{a.} If $F\cdot E=0$ for some 
exceptional curve $E$, consider the following cases.
\itemitemitemitem{i.} If $E\cdot L\ge 2$, then 
replace $F$ by $F-E$ and return to {\bf START}.
(Replacing $F$ by $F-E$ reduces 
$h^0(X, \C O_X(F))$ by 1, but by \ir{Elemma}
does not change the dimension of the kernel of $\mu_F$.)
\itemitemitemitem{ii.} If $E\cdot L=1$, 
then \ir{bnds}(e) gives the dimension of the kernel of
$\mu_F$: {\bf STOP}.
\itemitemitemitem{iii.} If $E\cdot L=0$, then
contracting $E$ gives
a birational morphism $\pi:X\to X'$, with respect to
which $\C O_X(F)=\pi^*\C O_{X'}(F')$ for some
(in fact canonically determined) $F'$, where
$F'$ is an effective divisor on $X'$ and fixed component free
and $\mu_{F}$ and $\mu_{F'}$ have kernels of the same dimension.
But $r$ has been reduced by 1, because $X'$ is a blow up of \pr2 at $r-1$
points. So replace $X$ by $X'$ and $F$ by $F'$,
and return to step 2.
\itemitemitem{b.} We thus reduce to the case 
that $F\ne0$ is effective, fixed component free and
has $F\cdot E>0$ for all exceptional curves $E$.
\itemitemitemitem{i.} If $F^2=0$, then 
\ir{pEncil} applies by \ir{introlem}, giving
$\hbox{dim}\,\hbox{ker}(\mu_F)$: {\bf STOP}.
\itemitemitemitem{ii.} If $F^2>0$, then $F$ is 
ample by \ir{introlem}: {\bf STOP}.

\irSubsection{The Main Result}{mnrs}
By \ir{amplesthm}, for $r\le 7$ general points of \pr2,
$\mu_F$ is surjective when $F$ is ample.
Thus, the algorithm above determines 
the rank of $\mu_F$ for an arbitrary $F$ on a blowing up of
\pr 2 at $r\le 7$ general points. But as mentioned above,
it is also desirable to have an explicit result 
in the case that $F$ is numerically effective. An
analysis of our algorithm
for numerically effective divisors
leads to an especially simple such result, \ir{mainthm}.

So say $r=7$. Denote $h^0(X, \C O_X(F+L))-3h^0(X, \C O_X(F))$
by $\lambda_F'$ and let $\lambda_F$ be the maximum of 0 and
$\lambda_F'$; note that $\lambda_F'=\hbox{dim}\,\hbox{cok}(\mu_F)-
\hbox{dim}\,\hbox{ker}(\mu_F)$ and that
$\mu_F$ has maximal rank
if and only if $\hbox{dim}\,\hbox{cok}(\mu_F)=\lambda_F$. 
Let $t_F$ be the number of
exceptional curves $E$ on $X$ with $E\cdot L=3$ such that $E\cdot F=0$.
It is well known (see [\refMa, \refNtwo])
that $[E]$ is the class of an exceptional curve with $E\cdot L=3$
if and only if $[E]$ is, up to permutation of the $E_i$, 
$[3L-2E_1-E_2-\cdots-E_7]$.
We denote these seven by $C_1=3L-2E_1-E_2-\cdots-E_7$,
$C_2=3L-E_1-2E_2-\cdots-E_7$, etc. We now have:

\prclm{Theorem}{mainthm}{Let $F$ be a numerically effective divisor
on the blow up $X$ of \pr2 at 7 general points,
$[L],[E_1],\cdots,[E_7]$ being the corresponding exceptional configuration. 
Then $\hbox{dim}\,\hbox{cok}(\mu_F)=\hbox{max }(t_F, \lambda_F)$, unless 
$[F]$ is, up to permutation of the $E_i$, either
0, $[B]$, $[B+C_4]$, $[B+C_4+C_5]$, 
$[B+C_4+C_5+C_6]$, $[B+C_4+C_5+C_6+C_7]$, $[G]$  or $[G+C_7]$, 
where $B=4L-2E_1-2E_2-2E_3-E_4-\cdots-E_7$ and
$G=5L-2E_1-\cdots-2E_6-E_7$, in which case $\mu_F$ is injective and
$\hbox{dim}\,\hbox{cok}(\mu_F)=\lambda_F$.} 

Although \ir{mainthm} does not explicitly address
the failure of $\mu_F$ to have maximal rank, 
it follows from \ir{mainthm} that $\mu_F$ fails to have
maximal rank if and only if $t_F>\lambda_F$ with $[F]$ not 
among the stated exceptions. (For an explicit example,
if $H$ is ample, then $\mu_F$ fails to have maximal rank for 
$F=H+(H\cdot C_i)C_i$: by \ir{amplesthm}, 
$\mu_H$ and hence $\mu_F$ are not injective
and thus $1=t_F=\hbox{dim}\,\hbox{cok}(\mu_F)$ by 
\ir{mainthm}.)

On the other hand,
as a corollary of \ir{mainthm} we see for a numerically effective 
$F\subset X$ that $\mu_F$ never fails by much to have maximal rank:
$\mu_F$ is never more than 7 short of maximal rank.

\irrnSection{Generalities}{gbkgrnd}
We first recall a useful exact sequence from [\refMu]. 
For sheaves \C F and \C G on $X$, we will
denote the kernel of the natural map $H^0(X,\C F)\otimes 
H^0(X, \C G)\to H^0(X,\C F\otimes \C G))$
by $\C R(\C F,\C G)$ and the cokernel by 
$\C S(\C F,\C G)$. When $\C F=\C O_X(F)$ and $\C G=\C O_X(G)$
for divisors $F$ and $G$ on $X$, we will, if it is convenient, 
just write $\C R(F,G)$ and $\C S(F,G)$.

\prclm{Proposition}{Mumford}{Let $C\subset X$ be a curve
on a smooth projective surface $X$, and
let $A$ and $B$ be divisors on $X$, so we have the exact sequence
$0\to \C O_X(A-C)\to \C O_X(A)\to\C O_{C}\otimes\C O_X(A)\to 0$. 
Then there is an exact sequence  $$0\to \C R(\C O_X(A-C), \C O_X(B))\to
\C R(\C O_X(A), \C O_X(B))\to\C R(\C O_{C}\otimes\C O_X(A),\C O_X(B)).$$
If the restriction homomorphisms 
$H^0(X, \C O_X(A))\to H^0(C,\C O_X(A)\otimes\C O_C)$
and $H^0(X, \C O_X(A+B))\to H^0(C,\C O_X(A+B)\otimes\C O_C)$
are surjective (for example, if $h^1(X, \C O_X(A-C))=0=
h^1(X,\C O_X(A+B-C))$), this extends to an exact sequence
$$\eqalign{0\to & \C R(\C O_X(A-C), \C O_X(B))\to
\C R(\C O_X(A), \C O_X(B))\to\C R(\C O_{C}\otimes\C O_X(A),\C O_X(B))\to\cr
&\C S(\C O_X(A-C), \C O_X(B))\to\C S(\C O_X(A), \C O_X(B))\to
\C S(\C O_{C}\otimes\C O_X(A),\C O_X(B))\to0.\cr}$$}

It will be helpful to have bounds on the dimensions of 
\C R and \C S.

\prclm{Proposition}{bnds}{Let $F$ be an effective divisor 
with $h^1(X,\C O_X(F))=0$ on the blowing up
$X$ of \pr2 at $r$ distinct points $p_1,\ldots,p_r$, let 
$[L], [E_1],\ldots,[E_r]$ be the corresponding exceptional configuration,
and assume that $F\cdot E_1\ge \cdots \ge F\cdot E_r$.
Let $d=F\cdot L$, $h=h^0(X,\C O_X(F))$, $l_i=h^0(X,\C O_X(F-(L-E_i)))$,
and $q_i=h^0(X,\C O_X(F-E_i))$.
\item{(a)} Then $\mu_F$ has maximal rank if 
and only if $\hbox{max}(0,2h-d-2)=\hbox{dim}\,\C R(F,L)$.
\item{(b)} If $h^1(X,\C O_X(F-(L-E_1)))=0=h^1(X,\C O_X(F-E_1))$, 
then $l_1+q_1=2h-d-2$.
\item{(c)} In any case, we have $\hbox{max}(0,2h-d-2)
\le \hbox{dim}\,\C R(F,L) \le l_1+q_1$.
\item{(d)} We also have $l_1+l_2
\le \hbox{dim}\,\C R(F,L) \le l_1+l_2+
h^0(X,\C O_X(F+(L-E_1-E_2)))-h^0(X,\C O_X(F))$.
\item{(e)} If $[L-E_1-E_2]$
is the class of an irreducible curve with
$F\cdot (L-E_1-E_2)=0$, then $\hbox{dim}\,\C R(F,L)=l_1+l_2$
and $\hbox{dim}\,\C S(F,L) = 
h^1(X,\C O_X(F-(L-E_1)))+h^1(X,\C O_X(F-(L-E_2)))$.}

\Prf \ir{bnds}(a,b,c) is just Corollary 4.2 of [\newfatpts].
Consider (d). If we choose coordinates $x$, $y$ and $z$ where $x$ and $y$
pass through $p_1$ and $y$ and $z$ through $p_2$,
then (from the proof of Lemma 4.1 of [\newfatpts])
$l_i$ is just the dimension of the kernel of
the restriction of $\mu_F$ to
$H^0(X, \C O_X(F))\otimes V_i\to H^0(X, \C O_X(F+L))$,
where $V_1$ is the vector space span of $x$ and $y$ in $H^0(X, \C O_X(L))$
and where $V_2$ is the vector space span of $z$ and $y$.
It is easy to see that these two kernels have only 0 in common;
this gives the lower bound of (d). For the upper bound,
it suffices to show $l_1+q_1\le l_1+l_2+
(h^0(X,\C O_X(F+(L-E_1-E_2)))-h^0(X,\C O_X(F)))$.
Since $[F-(L-E_2)]=[F-E_1-E]$, where $E$ is the effective divisor 
in the class $[L-E_1-E_2]$,
this follows from taking cohomology of
$0\to \C O_X(F-(L-E_2))\to \C O_X(F-E_1)\to 
\C O_E\otimes\C O_X(F-E_1)\to 0$,
using $\C O_E\otimes\C O_X(F-E_1)\cong \C O_E\otimes\C O_X(F+E)$ and
the fact that $h^1(X,\C O_X(F))=0$ implies that 
$h^0(X,\C O_X(F+(L-E_1-E_2)))-h^0(X,\C O_X(F))=
h^0(E,\C O_E\otimes\C O_X(F+E))$.

Finally consider (e); then $E$ is irreducible 
and hence a fixed component of
$|F+E|$, so (d) gives us $l_1+l_2=\hbox{dim}\,\C R(F,L)$.
From $\hbox{dim}\,\C S(F,L) = 
h^0(X, \C O_X(F+L) -3h^0(X, \C O_X(F)) +\hbox{dim}\,\C R(F,L)$,
we thus obtain
$\hbox{dim}\,\C S(F,L) = 
h^0(X, \C O_X(F+L)) -3h^0(X, \C O_X(F)) +
h^0(X,\C O_X(F-(L-E_1)))+h^0(X,\C O_X(F-(L-E_2)))$.
But $h^1(X, \C O_X(F))=0$ and hence $h^1(X, \C O_X(F+L))=0$
so Riemann--Roch gives $h^0(X, \C O_X(F+L))=h^0(X, \C O_X(F))
+F\cdot L+2$. Riemann--Roch also gives
$h^0(X,\C O_X(F-(L-E_i)))=h^0(X,\C O_X(F))+h^1(X,\C O_X(F-(L-E_i)))
-1-F\cdot (L-E_i)$ for $i=1,2$. 
Now substituting into our expression for $\hbox{dim}\,\C S(F,L)$
and simplifying (using $F\cdot L-F\cdot (L-E_1)-F\cdot(L-E_2)=-F\cdot E=0$)
gives the result.
\qed

\rem{Remark}{bndsrem} Note that the conclusion 
$\hbox{dim}\,\C R(F,L)=l_1+l_2$ of \ir{bnds}(e)
does not need the hypothesis that $h^1(X,\C O_X(F))=0$. The argument
that $l_1+l_2\le \hbox{dim}\,\C R(F,L)$ 
does not use $h^1(X,\C O_X(F))=0$,
and by Lemma 4.1 of [\newfatpts] neither does 
$\hbox{dim}\,\C R(F,L) \le l_1+q_1$.
Finally, with $E$ as in the proof of \ir{bnds}(e),  we have 
$h^0(E,\C O_E((F-E_1)\cdot E))=0$, so $l_2=q_1$ follows by 
taking cohomology of 
$0\to \C O_X(F-(L-E_2))\to \C O_X(F-E_1)\to \C O_E((F-E_1)\cdot E)\to 0$.
\qed

\prclm{Lemma}{Elemma}{Let $F\ne 0$ be an effective divisor 
on a smooth projective surface $X$,
and let $E$ be an exceptional curve with $F\cdot E=0$.
\item{(a)}  Say $|F|$ is fixed component free. Then 
$h^0(X,\C O_X(F-E))>0$, and
if $h^1(X,\C O_X(F))=0$, then  $h^1(X,\C O_X(F-E))=0$.
\item{(b)} Say $X$ is a blowing up of points of \pr2 and
$L$ is the total transform of a line.
If $E\cdot L\ge 2$, then the kernels
of $\mu_F$ and $\mu_{F-E}$ have the same dimension.}

\Prf (a) Since $F\cdot E=0$, we have an exact sequence
$0\to \C O_X(F-E)\to \C O_X(F)\to \C O_E\to 0$.
Since $|F|$ has no fixed components, $h^0(X, \C O_X(F))>1$
and $H^0(X, \C O_X(F))\to H^0(E, \C O_E)$ is surjective.
From the latter, 
our sequence is exact on global sections, so
our conclusions follow.

(b) Because $E\cdot L\ge 2$, it follows that
$h^0(X,\C O_X(L-E))=0$, but clearly 
$H^0(X,\C O_X(L-E))=\C R(\C O_E,\C O_X(L))$,
so $\C R(\C O_E,\C O_X(L))=0$.
Now apply \ir{Mumford} to the exact sequence 
in the proof of (a) to get an isomorphism
$\C R(F-E,L)\to\C R(F, L)$; i.e., the kernels
of $\mu_F$ and $\mu_{F-E}$ have the same dimension.\qed

\prclm{Lemma}{pEncil}{Let $X$ be a blowing up of distinct
points of \pr2 with corresponding exceptional configuration 
$[L], [E_1],\ldots,[E_r]$. Let 
$D\subset X$ be a smooth rational
curve with $D^2=0$ and let $m\ge 0$ be a nonnegative integer.
Then $\hbox{dim}\,\C R(mD, L)=m$ if $D\cdot L=1$ and 
$\C R(mD, L)=0$ if $D\cdot L>1$.}

\Prf If $L\cdot D>1$, then (as in the proof of \ir{Elemma}(b))
$0=H^0(X, \C O_X(L-D))=\C R(\C O_D,\C O_X(L))$. Applying 
\ir{Mumford} and induction on $s$
to $0\to \C O_X(sD)\to\C O_X((s+1)D)\to \C O_D\to 0$
gives $\C R(mD,L)=0$.

If $L\cdot D=1$, then  $[D]$ must be $[L-E_i]$ for some $i$.
By [\fatpts] (or directly), $\C S(mD,L)=0$, so
$\hbox{dim}\,\C R(mD,L)=3h^0(X,\C O_X(mD))-h^0(X,\C O_X(mD+L))=
3(m+1)-(2m+3)=m$. \qed

\irrnSection{Particularities}{sbkgrnd}
Now let $X$ be obtained by blowing up $r\le 8$ 
general points $p_1,\ldots,p_r$
of \pr2 and let $[L]$, $[E_1]$, $\ldots$, $[E_r]$ be the corresponding
exceptional configuration. We recall some facts for which we refer
to [\trans], [\antican], [\refHa] and [\refNtwo]. 

The exceptional configuration $[L]$, $[E_1]$, $\ldots$, $[E_r]$ is 
determined by and in turn determines a birational morphism $X\to \pr2$ 
with a factorization into monoidal transformations. 
Since $X$ can have more than one birational morphism to \pr2,
each of which typically factors in several ways,
$X$ can also have more than one exceptional configuration. For example,
if $\pi_1:X\to \pr2$ is the morphism determined by 
$[L]$, $[E_1]$, $\ldots$, $[E_r]$, and if 
$\pi_2:X\to \pr2$ is the morphism such that $\pi_2\pi_1^{-1}$ is
the quadratic Cremona transformation centered at $p_1$, $p_2$ and $p_3$
(i.e., $\pi_2\pi_1^{-1}$ is
the birational map from \pr2 to \pr2 given by the linear system of 
conics with base points at $p_1$, $p_2$ and $p_3$), then the 
exceptional configuration determined by $\pi_2$ (after an appropriate
factorization) is $[2L-E_1-E_2-E_3]$, $[L-E_2-E_3]$, $[L-E_1-E_3]$, 
$[L-E_1-E_2]$, $[E_4]$, $\ldots$, $[E_r]$.

Any two exceptional configurations are 
related by an element of the orthogonal
group on $\hbox{Cl}(X)$. Inside
the orthogonal group on $\hbox{Cl}(X)$, the subgroup $W$ 
generated by the reflections
$s_i$, $0\le i<r$, where $s_0(x)=
x+(x\cdot[L-E_1-E_2-E_3])[L-E_1-E_2-E_3]$
and $s_i(x)=x+(x\cdot[E_i-E_{i+1}])[E_i-E_{i+1}]$, is known as
the {\it Weyl\/} group. For $i>0$, the action
of $s_i$ on $[a_0L+a_1E_1\cdots +a_rE_r]$ is just to transpose the coefficients 
$a_i$ and $a_{i+1}$, while $s_0$ takes $[L]$, $[E_1]$, $\ldots$, $[E_r]$
to $[2L-E_1-E_2-E_3]$, $[L-E_2-E_3]$, $[L-E_1-E_3]$, $[L-E_1-E_2]$, 
$[E_4]$,$\ldots$, $[E_r]$. More generally, given 
any pair of exceptional configurations
there is an element of $W$ taking one to the other, and any $w\in W$ takes
$[L]$, $[E_1]$, $\ldots$, $[E_r]$ to another exceptional configuration.
This gives a bijection between exceptional configurations 
and elements of $W$.

If $[F_1]$ and $[F_2]$ are divisor classes
in the same orbit of $W$, then
$h^i(X, \C O_X(F_1))=h^i(X, \C O_X(F_2))$ holds for all $i$.
In addition, if $F$ is effective, then 
$h^2(X, \C O_X(F))=0$, while if $F$ is numerically effective,
then $h^1(X, \C O_X(F))=0$ and $|F|$ is nonempty and 
fixed component free.

Given any effective divisor $D$, we can write $[D]=[H]+[N]$, 
where $H$ is numerically effective and $N=-\sum(E\cdot D)E$,
where the sum is over all exceptional curves $E$ with $E\cdot D<0$;
note that the summands $E$ which appear in $N$ are disjoint.
Since, as noted above, $h^1(X,\C O_X(H))=0$, it is easy to
verify that $h^1(X,\C O_X(D))=0$ if and only if no summand 
in $N$ occurs with a coefficient of 2 or more (and hence
if and only if $D\cdot E\ge -1$ for every exceptional curve $E$).

For $8\ge r\ne2$, 
the classes of exceptional curves comprise one orbit, $W[E_r]$.
(If $r=2$, there are only three classes of exceptional curves,
$[L-E_1-E_2]$, $[E_1]$ and $[E_2]$, split between two $W$-orbits:
$\{[L-E_1-E_2]\}$ is one orbit, and $\{[E_1],[E_2]\}$ is 
the other.) For $r=7$, up to permutations of the $E_i$, 
the classes of the exceptional curves are just $[E_7]$,
$[L-E_1-E_2]$, $[2L-E_1-\cdots-E_5]$, and $[3L-2E_1-E_2-\cdots-E_7]$.

Also for $r=7$, the classes of 
numerically effective divisors are precisely
the $W$-orbits of nonnegative linear combinations of 
the classes of $L$, $L-E_1$, $2L-E_1-E_2$, 
$3L-E_1-E_2-E_3$, $\ldots$, $3L-E_1-\cdots-E_7$.
By excluding elements which can be obtained from others
(for example, exclude $[2L-E_1-E_2]$, since 
$[2L-E_1-E_2] =[L-E_1]+[L-E_2]$),
we can give a more efficient list of generators for the
cone of numerically effective divisor classes.
We thereby get the following list of divisors, whose classes
give a set of generators (complete up to permutation of the $E_i$) 
for the numerically effective cone:
\item{$G_1$}$=\div 1  -0  -0  -0  -0  -0  -0  -0 $,
\item{$G_2$}$=\div 2 -1 -1 -1  -0  -0  -0  -0 $,
\item{$G_3$}$=\div 3 -2 -1 -1 -1 -1  -0  -0 $,
\item{$G_4$}$=\div 4 -2 -2 -2 -1 -1 -1  -0 $,
\item{$G_5$}$=\div 4 -3 -1 -1 -1 -1 -1 -1 $,
\item{$G_6$}$=\div 5 -3 -2 -2 -2 -1 -1 -1 $,
\item{$G_7$}$=\div 5 -2 -2 -2 -2 -2 -2  -0 $,
\item{$G_8$}$=\div 6 -3 -3 -2 -2 -2 -2 -1 $,
\item{$G_9$}$=\div 7 -3 -3 -3 -3 -2 -2 -2 $,
\item{$G_{10}$}$=\div 8 -3 -3 -3 -3 -3 -3 -3 $,
\item{$G_{11}$}$=\div 1 -1  -0  -0  -0  -0  -0  -0 $,
\item{$G_{12}$}$=\div 2 -1 -1 -1 -1  -0  -0  -0 $,
\item{$G_{13}$}$=\div 3 -2 -1 -1 -1 -1 -1  -0 $,
\item{$G_{14}$}$=\div 4 -2 -2 -2 -1 -1 -1 -1 $,
\item{$G_{15}$}$=\div 5 -2 -2 -2 -2 -2 -2 -1 $,
\item{$G_{16}$}$=\div 3 -1 -1 -1 -1 -1 -1  -0 $,
\item{$G_{17}$}$=\div 4 -2 -2 -1 -1 -1 -1 -1 $,
\item{$G_{18}$}$=\div 5 -2 -2 -2 -2 -2 -1 -1 $, 
\item{$G_{19}$}$=\div 6 -3 -2 -2 -2 -2 -2 -2 $ and
\item{$G_{20}$}$=\div 3 -1 -1 -1 -1 -1 -1 -1 $. 

Since $[G_1]$ is clearly the class of a smooth rational
curve, so are $[G_2],\ldots, [G_{10}]$, since in fact they all
are in the same orbit of $W$. Likewise, $[G_{11}],\ldots, [G_{15}]$
is each the class of a smooth rational curve, and 
$[G_{16}],\ldots, [G_{20}]$
is each the class of a smooth elliptic curve.

It is also easy to check that each class $[G_i]$ is a sum of  
classes of exceptional curves and hence for $r=7$
the class of any effective divisor 
is a sum of classes of exceptional curves. It now follows for $r= 7$
(and in any case is well known) that $[3L-E_1-\cdots-E_7]=-K_X$ is ample
and hence so is any class of the form $[D]-K_X$, where $D$
is numerically effective. Conversely, for $r=7$
any ample class $[F]$ is of this form: as noted above, for some $w\in W$,
$w[F]$ is a nonnegative linear combination of the classes of 
the divisors $L$, $L-E_1$, $2L-E_1-E_2$, 
$3L-E_1-E_2-E_3$, $\ldots$, $3L-E_1-\cdots-E_7$.
But $w[F]\cdot [E_7]=[F]\cdot w^{-1}[E_7]>0$ since 
$w^{-1}[E_7]$ is the class of an exceptional curve and $F$ is ample,
so this linear combination involves $-K_X$
and hence is of the form $[D]-K_X$. I.e., 
$[F]=w^{-1}([D]-K_X)$, but $W$ preserves the numerically effective
cone, so in particular $w^{-1}[D]$ is numerically
effective. Finally, $w^{-1}(-K_X)=-K_X$ since $-K_X$
is stabilized by $W$, so $[F]$ has the required form.

The same argument works for $3\le r<7$; i.e., every
ample divisor class is $-K_X$ plus a numerically effective
class. The argument fails for $0\le r\le 2$ (for one thing, 
$-K_X$ is itself no longer needed as a generator
of the numerically effective cone if $0\le r\le 2$,
and, if $r<2$, the class of an effective divisor
need not be the sum of classes of exceptional curves).
However, an easy ad hoc argument shows that
the conclusion is still true for $r=2$, while for $r=1$
the ample divisor classes are $[dL-mE_1]$, where
$d>m$, and for $r=0$ they are $[dL]$, where
$d>0$.

\irrnSection{Application to 7 points}{mainres}
As an application of our results above, we will
prove \ir{mainthm}. To do so, we need some additional results.
We begin by considering ample divisors.

\prclm{Theorem}{amplesthm}{Let $F$ be an ample divisor 
on the blowing up $X$ of \pr2 at $t\le 7$ general points,
with $L$ the total transform of a line in \pr2.
Then $\C R(F,L)\ne 0$ and $\C S(F,L)=0$.}

\Prf Let $[L]$, $[E_1]$, $\ldots$, $[E_t]$ be the exceptional 
configuration corresponding to the $t$ 
points blown up to obtain $X$. After reindexing, we may assume that
$F\cdot E_1\ge F\cdot E_2\ge \cdots F\cdot E_t>0$.
If $t\le 5$, then, in fact, $\C S(F,L)=0$ for any numerically effective
$F$ by [\ir{fatpts}], while for $t\le2$, as follows from a
discussion above, every ample class $F$
is of the form $[L]$ plus a numerically effective
class. But $\C R(L,L)\ne 0$, so of course $\C R(F,L)\ne 0$, too.

So now we may assume $t\ge 3$.
Since $F$ is ample, as pointed out above we have $[F]=[D]-K_X$,
where $D$ is numerically effective. But $-K_X=[3L-E_1-\cdots-E_t]$, 
so $[F-E_1]=[D]+[3L-2E_1-\cdots-E_t]=[D+C_1]$. For $t<7$, $[C_1]$
is numerically effective and hence $0<h^0(X, \C O_X(F-E_1))=q_1$
and $0=h^1(X, \C O_X(F-E_1))$. If $t=7$, then
$[C_1]$ is the class of an exceptional curve, so 
$[F-E_1]$ is the class of an effective divisor, so
$0<h^0(X, \C O_X(F-E_1))=q_1$. 
Moreover, $E\cdot (D+C_1)\ge -1$ for every exceptional curve $E$,
so $F-E_1$ is regular (i.e., $h^1(X, \C O_X(F-E_1))=0$).

Since $h^1(X, \C O_X(F-E_1))=0$,
if we show $h^1(X,\C O_X(F-(L-E_1)))=0$, then by \ir{bnds} we will
know that $\mu_F$ has maximal rank and, using \ir{bnds} and
$q_1>0$ to see that $\C R(F,L)\ne 0$, that $\C S(F,L)$ must vanish.

From $[F]=[D]-K_X$ we obtain $[F-(L-E_1)]=[D]+[Q]$,
where $Q=2L-E_2-\cdots-E_t$. Arguing as for $h^1(X, \C O_X(F-E_1))$,
$h^1(X,\C O_X(F-(L-E_1)))$ also vanishes if $t<7$, so
we are reduced to the case that $t=7$.
Now, $[D]$ is a sum 
of classes $[U_i]$, where each divisor $U_i$ is, up to permutation
of the $E_i$, one of the divisors $G_j$ of \ir{sbkgrnd}. 
Recall each of the classes $[G_j]$ is the class of a smooth
curve, either rational or elliptic; by considering all permutations
of the $E_i$ for each $G_j$, we explicitly check
that $U_i\cdot(U_i+Q)\ge 2$ in each case that $[U_i]$ is 
the class of an elliptic curve and 
$U_i\cdot(U_i+Q)\ge -1$ in each case that $[U_i]$ is the class
of a rational curve, 
unless $U_i=\div 5 -1 -2 -2 -2 -2 -2 -2 $, in which case
$U_i\cdot(U_i+Q)= -2$. Thus, letting $A_i$
be a smooth curve with $[A_i]=[U_i]$, we have 
$h^1(A_i, \C O_{A_i}(U_i+Q))=0$
unless $U_i=\div 5 -1 -2 -2 -2 -2 -2 -2 $. Moreover,
$(\div 5 -1 -2 -2 -2 -2 -2 -2 )\cdot U_i>0$ for all $i$ 
with $U_i\ne \div 5 -1 -2 -2 -2 -2 -2 -2 $. Thus, unless each
$U_i$ is $\div 5 -1 -2 -2 -2 -2 -2 -2 $, we may assume that $U_1$ is not
$\div 5 -1 -2 -2 -2 -2 -2 -2 $, and then from $h^1(X,\C O_X(Q))=0$
it follows inductively by taking cohomology of 
$$0\to \C O_X(Q+U_1+\cdots+U_{i-1})\to \C O_X(Q+U_1+\cdots+U_i)\to
\C O_{A_i}(Q+U_1+\cdots+U_i)\to 0$$
that $h^1(X,\C O_X(D+Q))=0$, as desired.

There remains the case that $F=m(\div 5 -1 -2 -2 -2 -2 -2 -2 )-K_X$,
for $m>0$. But our assumption that 
$F\cdot E_1\ge F\cdot E_2\ge \cdots F\cdot E_7>0$ rules out this case.
\qed

\prclm{Lemma}{crit}{Let $X$ be a blowing up of \pr2 at 7 general points
$p_1,\ldots,p_7$, with $[L]$, $[E_1],\ldots,[E_7]$ 
the corresponding exceptional
configuration. Let $J_i$, $i=1,2$, be smooth curves whose classes are
$[L-E_i]$. Let $0\ne [F]$ 
be numerically effective with $F\cdot (L-E_1-E_2)=0$ and 
$F\cdot E_1\ge \cdots \ge F\cdot E_7$. 
Then $\mu_F$ fails to have maximal rank if and only if 
$h^0(X,\C O_X(F-J_1))>0$ and $h^1(X,\C O_X(F-J_2))>0$.}

\Prf By \ir{bnds}(e), $l_1>0$ implies that $\mu_F$ is not injective, while
$h^1(X,\C O_X(F-J_2))>0$ implies that $\mu_F$ is not surjective.

Conversely, by \ir{bnds}(e), if $\mu_F$ is neither surjective nor injective,
then $l_1+l_2>0$ and $h^1(X,\C O_X(F-J_1))+h^1(X,\C O_X(F-J_2))>0$,
so it suffices to check that $l_2>0$ implies $l_1>0$, and that
$l_1>0$ and $h^1(X,\C O_X(F-J_1))>0$ 
together imply $h^1(X,\C O_X(F-J_2))>0$.

Suppose $l_2>0$. Thus $[F-J_2]$ is a sum of classes of exceptional 
curves $T_i$, and, since $F\cdot E_1\ge F\cdot E_2$ and hence
$(F-J_2)\cdot E_1> (F-J_2)\cdot E_2$, some summand
has $[T_i]\cdot (E_1-E_2) >0$, hence by Riemann--Roch
and duality $h^0(X,\C O_X(T_i+(E_1-E_2)))>0$. Thus 
$l_1=h^0(X,\C O_X(F-J_2+(E_1-E_2)))>0$, as claimed.

Now assume $l_1>0$ and $h^1(X,\C O_X(F-J_1))>0$. 
If $F\cdot E_1= F\cdot E_2$, then $F-J_1$ and $F-J_2$ are the same,
up to permutation of the $E_i$, hence
in the same orbit of the Weyl group, so 
$h^1(X,\C O_X(F-J_1))=h^1(X,\C O_X(F-J_2))$.
So suppose that $F\cdot E_1>F\cdot E_2$, and hence that
$(F-J_1)\cdot E_1\ge \cdots \ge (F-J_1)\cdot E_7\ge 0$.
Since $[F-J_1]$ 
has an effective representative, $h^1(X,\C O_X(F-J_1))>0$ implies
that there is an exceptional curve $E$ with $(F-J_1)\cdot E\le -2$.
Clearly, this $E$ is not among the $E_i$, so we may assume that
$[E]$ is either $[L-E_1-E_2]$, 
$[2L-E_1-\cdots-E_5]$ or $[3L-2E_1-E_2-\cdots-E_7]$
(since up to permutation of the $E_i$, the class of 
every exceptional curve is one of these, 
and these are the permutations minimizing the intersection with $F-J_1$).
But whichever of these is $E$, we have 
$(F-J_2)\cdot E=(F-J_1-(E_1-E_2))\cdot E
\le (F-J_1)\cdot E\le -2$, so from
$0\to \C O_X(F-J_2-E)\to \C O_X(F-J_2)\to \C O_E((F-J_2)\cdot E)\to 0$,
it suffices to check that $h^2(X,\C O_X(F-J_2-E))=0$ to
obtain that $h^1(X,\C O_X(F-J_2))>0$, as required.
But $F\cdot L\ge 1$ (since $[F]$ is nontrivial and numerically effective),
so $(K_X-[F-J_2-E])\cdot L<0$ (so $0=h^0(X,\C O_X(K_X-(F-J_2-E)))
=h^2(X,\C O_X(F-J_2-E))$, since $L$ is numerically effective) 
unless $E\cdot L=3$ and $F\cdot L=1$. In this latter case 
$[E]=[3L-2E_1-E_2-\cdots-E_7]$ and $[F]=[L-E_1]$. 
Then $K_X-[F-J_2-E]=[-E_2]$,
which again is not the class of an effective divisor, so
again $h^2(X,\C O_X(F-J_2-E))=0$ by duality. \qed

\prclm{Lemma}{mxrkEperp}{Let $X$ be as in \ir{crit},
let $E$ be an exceptional curve with $E\cdot L=1$, 
and let $F$ be numerically effective
such that $F\cdot E=0$, but $F\cdot C>0$ for every exceptional
curve $C$ with $C\cdot L\ne1$. 
Then $\C S(F, L)=0$ but $\C R(F, L)\ne 0$.}

\Prf As usual, we may assume that 
$F\cdot E_1\ge \cdots \ge F\cdot E_7$, and thus we may assume
$E$ is the exceptional curve whose class is
$[L-E_1-E_2]$.

First say that $F\cdot C>0$ for every exceptional curve $C\ne E$.
Choose an element $w$ of the Weyl group $W$
such that $w[F]$ is a sum of
nonnegative multiples of the classes of $L$, $L-E_1$, $2L-E_1-E_2$, 
$3L-E_1-E_2-E_3$, $\ldots$, $3L-E_1-\cdots-E_7$. 
Note that this sum cannot involve $-K_X=[3L-E_1-\cdots-E_7]$.
(If it did, then $w[F]=[D]-K_X$ for some numerically effective $D$,
but $-K_X$ is ample and hence so would be $w[F]$ and thus $[F]$,
contradicting $F\cdot E=0$.) It follows that $w[F]\cdot E_7=0$
and hence that $w[E]=[E_7]$. Since $F\cdot C>0$ for every
exceptional curve $C\ne E$, we have $w[F]\cdot E_6>0$, hence the class of
$H=3L-E_1-\cdots-E_6$ appears in the sum. Thus $[F]-w^{-1}[H]$ is 
numerically effective, so $F\cdot E=0$ implies $w^{-1}[H]\cdot E=0$.
But looking over the $W$-orbit of $[H]$ shows 
it has only one element perpendicular to $E$; i.e.,  we must have 
$w^{-1}[H]=[\div 4 -2 -2 -1 -1 -1 -1 -1 ]$. Thus 
$[(F-w^{-1}[H])+(w^{-1}[H]-(L-E_1))]= [D+C_2]$ 
for some numerically effective $D$.
Since  $[C_2]=[\div 3 -1 -2 -1 -1 -1 -1 -1 ]$ is 
the class of an exceptional curve,
we see that $l_1=h^0(X, \C O_X(F-(L-E_1)))>0$ 
and $h^1(X, \C O_X(F-(L-E_1)))=0$;
similarly, $h^1(X, \C O_X(F-(L-E_2)))=0$. By \ir{bnds}(e), 
$\C R(F, L)\ne 0$ and $\C S(F, L)=0$.

Now suppose that $F\cdot C=0$ for some exceptional curve $C\ne E$.
If we denote $L-E_i-E_j$ by $C_{ij}$, then
by hypothesis $[C]=[C_{ij}]$ for some $i$ and $j$, and, since 
$F\cdot E_1\ge \cdots \ge F\cdot E_7$,
either $F\cdot E_1= F\cdot E_2$ and thus $[F]$
is of the form $[2a_1L-a_1(E_1+\cdots +E_i)-b_{i+1}E_{i+1}-\cdots-b_7E_7]$
where $a_1>b_{i+1}\ge \cdots\ge b_7>0$ and $i\ge 3$, or 
$F\cdot E_1> F\cdot E_2$ and thus $[F]$
is of the form $[(a_1+a_2)L-a_1E_1-a_2(E_2+\cdots +E_i)-
b_{i+1}E_{i+1}-\cdots-b_7E_7$
where $a_1>a_2>b_{i+1}\ge \cdots\ge b_7>0$ and $i\ge 3$.

For the former, $i=3$, since otherwise $F\cdot (2L-E_1-\cdots -E_5)\le0$,
so the classes of the only exceptional curves that 
$F$ is perpendicular to are
$[C_{12}]$, $[C_{13}]$, and $[C_{23}]$. As above, 
$w\{C_{12}, C_{13}, C_{23}\}
=\{E_5,E_6,E_7\}$ and $w[F]$ is a nonnegative sum of the classes of 
$L$, $L-E_1$, $2L-E_1-E_2$, 
$3L-E_1-E_2-E_3$, and $3L-E_1-E_2-E_3-E_4$,
for some $w\in W$, and this sum involves
$H=3L-E_1-E_2-E_3-E_4$. Thus $w^{-1}[H]$ is perpendicular
to each of $[C_{12}]$, $[C_{13}]$, and $[C_{23}]$, but
by examining the $W$-orbit of $H$, we see there is only one element of
$W[H]$ perpendicular to each of $[C_{12}]$, $[C_{13}]$, and $[C_{23}]$;
i.e., $w^{-1}[H]=[\div 6 -3 -3 -3 -1 -1 -1 -1 ]$. But 
$w^{-1}[H]-[L-E_1] = [C_{23}] + 
[\div 4 -2 -2 -2 -1 -1 -1 -1 ]$ and $w^{-1}[H]-[L-E_2] = [C_{13}] + 
[\div 4 -2 -2 -2 -1 -1 -1 -1 ]$; since $\div 4 -2 -2 -2 -1 -1 -1 -1 $
is numerically effective, we conclude
that $l_1=h^0(X, \C O_X(F-(L-E_1)))>0$
and $h^1(X, \C O_X(F-(L-E_1)))=0=h^1(X, \C O_X(F-(L-E_2)))$ and hence
$\C R(F, L)\ne 0$ and $\C S(F, L)=0$ by \ir{bnds}(e).

For the latter, $F$ is perpendicular to $C_{1j}$
for all $2\le j\le i$, where, we recall, $i\ge 3$. Reasoning
as above, for some $w\in W$, $[F]$ is a sum of a numerically effective
class $[D]$ and $w^{-1}[M_i]$, where $M_i=3L-E_1-\cdots -E_{8-i}$
for $3\le i\le 6$ and $M_7=2L-E_1$, 
and where $w^{-1}[M_i]$ is perpendicular to each $C_{1j}$
but to no other exceptional curves. 
As above, there is in each case a unique possibility
for $w^{-1}[M_i]$:
$w^{-1}[M_3]=[\div 5 -3 -2 -2 -1 -1 -1 -1 ]$;
$w^{-1}[M_4]=[\div 6 -4 -2 -2 -2 -1 -1 -1 ]$;
$w^{-1}[M_5]=[\div 7 -5 -2 -2 -2 -2 -1 -1 ]$;
$w^{-1}[M_6]=[\div 8 -6 -2 -2 -2 -2 -2 -1 ]$; and
$w^{-1}[M_7]=[\div 5 -4 -1 -1 -1 -1 -1 -1 ]$.

In each of the cases $3\le i\le 6$ 
one checks as above that $w^{-1}[M_i]-[(L-E_1)]$ and hence $[F-(L-E_1)]$ 
are classes of effective divisors, and similarly
that $[F-(L-E_2)]$ is the class of an effective divisor 
with $[F-(L-E_2)]\cdot C\ge -1$ for every exceptional curve $C$.
This implies that $l_1>0$ and $h^1(X,\C O_X(F-(L-E_2)))=0$, 
as required. 

We are left with the case $[H]=w^{-1}[M_7]$.
We note that $h^0(X,\C O_X(H-(L-E_1)))>0$
and $h^1(X,\C O_X(H-(L-E_2)))=0$, but $h^0(X,\C O_X(H-(L-E_2)))=0$. 
By Riemann--Roch, $h^0(X,\C O_X(H+D-(L-E_2)))\ge 
h^0(X,\C O_X(D)) -1+D\cdot (H-(L-E_2))$. By checking each of the 
generators $[G_i]$ of the numerically effective cone (including
those obtained by permutations of the $E_i$), 
we see that $h^0(X,\C O_X(H+D-(L-E_2)))$
is positive unless $[D]$ is a nonnegative multiple
of $[\div 3 -2 -0  -1 -1 -1 -1 -1 ]$, in which case
$(H+D-(L-E_2))\cdot (\div 3 -2 -0  -1 -1 -1 -1 -1 )=-1$,
so numerical effectivity of $[\div 3 -2 -0  -1 -1 -1 -1 -1 ]$
implies $h^0(X,\C O_X(H+D-(L-E_2)))=0$, and now Riemann--Roch 
gives $h^1(X,\C O_X(H+D-(L-E_2)))=0$, as required.
If $[D]$ is not a multiple of $[\div 3 -2 -0  -1 -1 -1 -1 -1 ]$,
then $h^0(X,\C O_X(H+D-(L-E_2)))>0$, but then
$(H-(L-E_2))\cdot B\ge -1$ for every exceptional curve $B$
and hence the same is true for $D+H-(L-E_2)$ so again
$h^1(X,\C O_X(H+D-(L-E_2)))=0$.
\qed

\prclm{Lemma}{fccase}{Let $X$ be a blowing up of
\pr2 at 7 general points, $[L],[E_1],\cdots,[E_7]$
the corresponding exceptional configuration.
Let $[F]$ be a nontrivial numerically effective class and
let $E$ be an exceptional curve with $E\cdot F=0$.
If $C$ is a reduced irreducible curve occurring as a
fixed component of $|F-E|$, then $C$ is an exceptional curve,
$F^2=0$ and $[F]=m[E+C]$ for some $m>0$. In addition, if
$L\cdot (E+C)>1$, then $\C R(F,L)=0$.}

\Prf Suppose $C$ is a fixed component of $|F-E|$
(recall by \ir{Elemma}(a) that $|F-E|$ is nonempty).
Any integral curve $C$ is either numerically effective
or has $C^2<0$. But on a 7 point blow up, the former are never fixed
and the latter are exceptional; thus $C$ must be an exceptional curve.

Since $C$ is in the base locus of $|F-E|$, we can write
$[F-E]=[H]+[N]$, where $N$ and $H$ are the fixed and free parts,
respectively, of $|F-E|$ and $C$ is a component of $N$,
hence $C\cdot (F-E)=C\cdot N<0$, but $F$ is numerically effective
so $C\cdot E>0$. On the other hand, $E\cdot (H+N+E)=E\cdot F=0$, so
$E\cdot (H+N)=1$. Now, $E\cdot C>0$ implies that $|E+C|$
is positive dimensional, hence cannot be contained in $N$.
Of course, $C$ is in $N$, so $E$ cannot be. Thus $E\cdot N>0$,
so $E\cdot (H+N)=1$ tells us that $E\cdot H=0$ and $E\cdot N=1$.
Therefore, $E$ is perpendicular to components of $N$ other than
$C$ while $E\cdot C=1$ (which means that $|E+C|$ is a pencil). 
Since this would mean components of
$N$ other than $C$ would meet $[H+N+E]=[F]$ negatively, there can be
no other components and we see that $N=C$. Thus $H$ is perpendicular
to both $C$ and $E$, and therefore $|H|$ is composed with the pencil
$|E+C|$; i.e., $[H]$ is a multiple of $[E+C]$, so
$[F]=m[E+C]$ for some $m>0$.

Now let $L\cdot (E+C)>1$; then apply \ir{pEncil} with
$D$ a general element of $|E+C|$ to obtain $\C R(F,L)=0$. \qed

We now give the proof of \ir{mainthm}.

\Prf By the algorithm discussed in \ir{intro},
one can explicitly check that $t_F>\lambda_F=\hbox{dim}\,\C S(F,L)$ 
and $\C R(F,L)=0$ for each 
exception $F$ listed in the statement of the theorem.

We now show that otherwise $\hbox{dim}\,\C S(F,L)$ is the maximum of
$t_F$ and $\lambda_F$.
So let $F$ be a nontrivial numerically effective divisor.

It may be that $F\cdot E'=0$ for some exceptional curve $E'$ 
with $E'\cdot L\ge 2$. By \ir{Elemma}(a), $|F'|$ is
nonempty for $F'=F-E'$. We continue in this way, subtracting
off exceptional curves meeting $L$ at least twice,
to obtain a sequence $F=F', F'',\ldots, F^{(j)},\dots$, as long as
$F^{(j)}$ continues to be perpendicular to some such exceptional curve
$E^{(j+1)}$ and as long as $F^{(j+1)}=F^{(j)}-E^{(j+1)}$ 
has a fixed component free
linear system. Eventually, however, say for $j=t$, either 
$F^{(t)}$ is numerically 
effective with $F^{(t)}\cdot E>0$ for every exceptional curve $E$ 
with $E\cdot L\ge 2$, or $F^{(t)}$ is 
effective but $|F^{(t)}|$ has a fixed component.
For convenience, we write $Q_t$ for $E'$, $Q_{t-1}$ for $E''$,
etc., and also $F_0$ for $F^{(t)}$ so $[F^{(j)}]=
[F_{t-j}]$, where $F_{t-j}=F_0+Q_1+\cdots+Q_{t-j}$ for $0\le j \le t$.
Thus $F_i$ is numerically 
effective with $F_i\cdot Q_i=0$ for $i>0$, and $F_0$ is either numerically 
effective with $F\cdot E>0$ for every exceptional curve $E$ 
with $E\cdot L\ge 2$, or $F_0$ is 
effective but $|F_0|$ has a fixed component.

Since $0=H^0(X, \C O_X(L-Q_i))=\C R(\C O_{Q_i},\C O_X(L))$, we have
$\hbox{dim}\,\C S(\C O_{Q_i},\C O_X(L))=
h^0(Q_i,\C O_{Q_i}(Q_i\cdot L))-3=Q_i\cdot L-2$, so
applying \ir{Mumford} and induction on $i$ 
to $0\to \C O_X(F_{i-1})\to \C O_X(F_i)\to \C O_{Q_i}\to 0$
we see $\hbox{dim}\,\C S(F,L)=
\hbox{dim}\,\C S(F_0,L)+(Q_1+\cdots+Q_t)\cdot L-2t$. (Note that
$(Q_1+\cdots+Q_t)\cdot L-2t$ is just the number
of summands $Q_i$ with $Q_i\cdot L=3$.)

Consider first the case that $F_0$ is numerically effective.
If $F_0\cdot E_i=0$ for some $i$, then we can regard $F_0$ as a 
divisor on a blowing up of \pr2 at 6 points. By
[\refFi], $\C S(H,L)=0$ for all numerically effective divisors $H$
on a blowing up of \pr2 at 6 general points $p_1,\ldots,p_6$ unless $H$ is
$5L-2E_1-\cdots-2E_6$ or a multiple of $3L-2E_{i_1}-E_{i_2}-\cdots-E_{i_6}$.
But $F_0$ cannot be any of these since they are perpendicular
to exceptional curves meeting $L$ at least twice. Thus $\C S(F_0,L)=0$
if $F_0\cdot E_i=0$ for some $i$. Otherwise, $F_0\cdot E>0$ for every 
exceptional curve $E$ with $E\cdot L\ne 1$, and either $F_0$ is ample
(whence $\C S(F_0,L)=0$ by \ir{amplesthm}) or $F_0\cdot E=0$ for some 
exceptional curve $E$ with $E\cdot L=1$
(whence $\C S(F_0,L)=0$ by \ir{mxrkEperp}). Either way, we have 
$F_0\cdot Q_i>0$ for all $i$ and 
$\hbox{dim}\,\C S(F,L)=(Q_1+\cdots+Q_t)\cdot L-2t$. Since
$F_0\cdot Q_i>0$ but $F_i\cdot Q_i=0$, it follows inductively
that $Q_i\cdot Q_j=0$ for all $i\ne j$. It now follows easily that
$t_F=(Q_1+\cdots+Q_t)\cdot L-2t$; since $\lambda_F\le \hbox{dim}\,\C S(F,L)$
is always true, we have $\hbox{dim}\,\C S(F,L)=t_F=\hbox{max}(t_F,\lambda_F)$,
as claimed.

Now consider the case that $|F_0|$ has a fixed component. 
By \ir{fccase}, $[F_0+Q_1]$ is $m[H]$, where $H^2=0$
and $|H|$ is a pencil. Thus, up to indexation, $[H]$ is among
$[G_{11}]=[L-E_1]$, $[G_{12}]=[2L-E_1-\cdots-E_4]$, 
$[G_{13}]=[3L-2E_1-E_2-\cdots-E_6]$,
$[G_{14}]=[\div 4 -2 -2 -2 -1 -1 -1 -1 ]$, or
$[G_{15}]=[\div 5 -2 -2 -2 -2 -2 -2 -1 ]$, but 
$Q_1\cdot L\ge 2$ rules out $[L-E_1]$. If $t>1$, then we have 
$(F_0+Q_1)\cdot Q_2=1$, hence $m=1$. Thus $[F]$ is either
$m[H]$ or $[H+Q_2+\cdots+Q_t]$. In the former case
$\C R(F,L)=0$ by \ir{fccase}, hence 
$\hbox{dim}\,\C S(F,L)=\lambda_F$, and we explicitly check
that $t_F\le \lambda_F$ unless $m=1$ and $[H]$
is either $[\div 4 -2 -2 -2 -1 -1 -1 -1 ]$ or
$[\div 5 -2 -2 -2 -2 -2 -2 -1 ]$.

This verifies the statement of \ir{mainthm}
unless $[F]$ is of the form $[H+Q_2+\cdots+Q_t]$, as above,
where $[H]$ is one of $[2L-E_1-\cdots-E_4]$, $[3L-2E_1-E_2-\cdots-E_6]$,
$[\div 4 -2 -2 -2 -1 -1 -1 -1 ]$, or
$[\div 5 -2 -2 -2 -2 -2 -2 -1 ]$. We consider each possibility
for $[H]$ in turn. In each case we have $0=\C R(H,L)=\C R(F,L)$,
so $\hbox{dim}\,\C S(F,L)=\lambda_F$, and it is enough to check
that $t_F\le \hbox{dim}\,\C S(H,L) + \rho$, where
$\rho=(Q_2+\cdots+Q_t)\cdot L-2t+2$, when $F$ is not
one of the stated exceptions. 

First consider $[H]=[2L-E_1-\cdots-E_4]$. Since any exceptional $E$
with $E\cdot L=3$ has $E\cdot H>0$, we can have $E\cdot F=0$ only
if $E$ is among the $Q_i$. Thus $t_F\le \rho$, settling this case.

Now let $[H]=[3L-2E_1-E_2\cdots-E_6]$. Any exceptional $E$
with $E\cdot L=3$ and $E\cdot F=0$ must be among the $Q_i$, or
must have $E\cdot H=0$ (and hence $[E]=[3L-2E_1-E_2\cdots-E_7]$). 
Thus $t_F\le \rho+1$, but $\hbox{dim}\,\C S(H,L) =1$ in this case,
so this case is also settled.

Now suppose $[H]=[\div 4 -2 -2 -2 -1 -1 -1 -1 ]$. Let us say
that a divisor $B$ is {\it cubic\/} if $B\cdot L=3$ and
{\it conic\/} if $B\cdot L=2$. Now argue
as in the preceding paragraph. This time
there are exactly three cubic exceptionals perpendicular to $H$
(in fact, their classes $[C_i]$ are
exactly $-K_X-[E_i]$, $1\le i\le 3$),
so we see $t_F\le \rho+3$. 
Since now $\hbox{dim}\,\C S(H,L) =2$,
this case is settled unless $t_F=\rho+3$
and hence the classes of $C_1,C_2,C_3$ and of all of the cubics
among $Q_i$, $i\ge2$, are distinct and 
perpendicular to $F$ (else we certainly would have
$t_F<\rho+3$), in which case these $t_F$ cubics are also  
perpendicular to each conic $Q_i$, $i\ge 2$.
But the class of a conic exceptional curve perpendicular to
$C_1,C_2,\hbox{ and }C_3$  must be of the form 
$[2L-E_1-E_2-E_3-E_{i_1}-E_{i_2}]$ and is therefore
perpendicular to $H$, and hence the conic
$Q_j$ with least $j\ge 2$ must meet one of the cubics
occurring among the $Q_i$, $i\ge 2$. To avoid this contradiction
we conclude there is no conic $Q_i$, $i\ge 2$.
Thus $[F]$ must be $[H]$ plus any of the four classes of cubic
exceptionals not perpendicular
to $H$, giving only the exceptions
$[H]$, $[H+C_4]$, $[H+C_4+C_5]$, 
$[H+C_4+C_5+C_6]$, $[H+C_4+C_5+C_6+C_7]$ listed in the statement
of \ir{mainthm}.

Finally, we have $[H]=[\div 5 -2 -2 -2 -2 -2 -2 -1 ]$.
Here we have $\lambda_F=\hbox{dim}\,\C S(F,L)=\rho+3$,
so for $t_F>\lambda_F$ we need $t_F\ge 4$. 
Let us look at all numerically effective classes perpendicular to
at least four (say to $[C_4],\ldots, [C_7]$) cubic exceptionals. 
From the generators $[G_i]$
of the numerically effective cone given
in \ir{sbkgrnd} it is easy to verify that any numerically effective
class perpendicular to each of $C_4,\ldots, C_7$ is a nonnegative
sum of the classes of $A=\div 7 -2 -2 -2 -3 -3 -3 -3 $,
$D_i=\div 5 -2 -2 -2 -2 -2 -2 -2 +E_i$, $1\le i\le3$, and $B=
\div 8 -3 -3 -3 -3 -3 -3 -3 $.
We will show that any sum $F$ of these 5 divisors with $t_F>\lambda_F$
and $\C R(F,L)=0$ must be among the list of exceptions given
in the statement of \ir{mainthm}.

First note that the class of
each of these five divisors is on the list of exceptions.
Now let $F$ be $A$ plus any one of $A$, $D_1$, $D_2$, $D_3$ and $B$.
In each case we check that $h^0(X, \C O_X(F-C_4-C_5-C_6-C_7))>0$
and $h^1(X, \C O_X(F-C_4-C_5-C_6-C_7))=0$. Thus the same is true for
any sum $F$ of two or more
of the divisors $A$, $D_1$, $D_2$, $D_3$ and $B$ such that at 
least one summand is $A$. Applying \ir{Mumford} to
$0\to \C O_X(F-Y) \to \C O_X(F) \to \C O_Y\to 0$,
where $Y$ is the disjoint union of the $t_F=4$ cubic exceptional curves 
perpendicular to $F$ (so $[Y]=[C_4]+\cdots+[C_7]$), we conclude that
$\hbox{dim}\,\C S(F,L)=\hbox{dim}\,\C S(F-Y,L)+t_F$. 
Thus, whenever we have $\C R(F,L)=0$ for such an $F$, we also have
$\lambda_F=\hbox{dim}\,\C S(F, L)\ge t_F$. 

Now consider the case that $F$ is a sum with three or more summands
taken from $D_1$, $D_2$, $D_3$ and $B$. In every such case of
3 summands (and hence also for more than 3 summands) except for pure
multiples of some $D_i$ (which were treated above),
$[F]$ is, as in the preceding paragraph, 
the class of the sum of the cubic exceptionals perpendicular to $F$
plus an effective regular divisor, and therefore as above
$\C R(F,L)=0$ implies $t_F\le \lambda_F$. 

We are left to consider the case that $F$ is a sum of any two
of $D_1$, $D_2$, $D_3$ and $B$ except pure multiples of
some $D_i$: if $F=D_j+D_i$, $j\ne i$, then $t_F=5$ and $\lambda_F=4$;
if $F=B+D_i$, then $t_F=6$  and $\lambda_F=5$; 
and if $F=2B$, then $t_F=7$  and $\lambda_F=6$. But in each of these
cases, $F$ is one of the exceptions explicitly given in
the statement of the theorem. \qed

Although \ir{mainthm} is well-suited for computational applications;
our final result is more conceptually satisfying.

\prclm{Corollary}{lastcor}{Let $F$ be a numerically effective divisor on
$X$, with $X$ as in \ir{mainthm}. Let $D=C_{i_1}+\cdots+C_{i_{t_F}}$ 
be the sum of the cubic exceptional curves 
$C_{i_j}\in\{C_1,\ldots,C_7\}$ perpendicular to $F$.
If $\mu_F$ fails to have maximal rank, then $\hbox{dim}\,\C S(F,L)=t_F$.
Moreover, $\mu_F$ fails to have maximal rank if and only if:
$t_F>0$, $F-D$ is numerically effective, and $\lambda_{F-D}'<0$.}

\Prf If $\mu_F$ fails to have maximal rank, then 
$\hbox{dim}\,\C S(F,L)=t_F$ follows by \ir{mainthm}.
We now consider the second claim.

If $t_F=0$, then $\hbox{dim}\,\C S(F,L)=\lambda_F$ by \ir{mainthm},
so $\mu_F$ has maximal rank.

Consider the case that $t_F>0$ but 
$F-D$ is not numerically effective. 
Since $F-D$ is not numerically effective, successively
subtracting $C_{i_1},C_{i_2},\ldots$ from $F$, 
we eventually obtain by \ir{Elemma} (as in the proof of \ir{mainthm})
a divisor $F_0$ whose class is the class of an effective but not
numerically effective divisor. Now, by the proof of \ir{mainthm},
$S(F,L)=\lambda_F$ and thus $\mu_F$ has maximal rank.

Finally, say $t_F>0$ and $F-D$ is numerically effective.
Then $t_{F-D}=0$, so as we saw above, $\mu_{F-D}$
has maximal rank; in particular, $\mu_{F-D}$
has a nontrivial kernel if and only if $\lambda_{F-D}'<0$.
Now apply \ir{Mumford} to the exact sequence 
$0\to \C O_X(F-D)\to \C O_X(F)\to \C D\to 0$.
Since $\C R(\C O_D,\C O_X(L))=0$, we see that  
$\hbox{dim}\,\C R(F-D,L)=\hbox{dim}\,\C R(F,L)$
and $\hbox{dim}\,\C S(F,L)=\hbox{dim}\,\C S(F-D,L)
+ \hbox{dim}\,\C S(\C O_D,\C O_X(L))=\hbox{dim}\,\C S(F-D,L)+t_F\ge t_F>0$.
Thus $\mu_F$ fails to have maximal rank if and only
if $\mu_{F-D}$ fails to be injective, which
we noted above holds if and only if $\lambda_{F-D}'<0$.
\qed

\vskip-\baselineskip
\vskip-\baselineskip
\References

\bibitem{\refCat} Catalisano, M.\ V. {\it ``Fat'' points on a conic}, 
Comm.\ Alg.\  19(8) (1991), 2153--2168.

\bibitem{\refFi} Fitchett, S. {\it Doctoral dissertation}, 
University of Nebraska-Lincoln, 1997.

\bibitem{\refGO} Geramita, A. V. and Orrechia, F. {\it Minimally
generating ideals defining certain tangent cones}, J.\ Alg. 78
(1982), 36--57.

\bibitem{\refGGR} Geramita, A.\ V., Gregory, D.\ and Roberts, L.
{\it Minimal ideals and points in projective space},
J.\ Pure and Appl.\ Alg. 40 (1986), 33--62.

\bibitem{\trans} Harbourne, B. {\it Complete linear 
systems on rational surfaces}, 
Trans.\ A.\ M.\ S.\ 289 (1985), 213--226. 

\bibitem{\vanc} \manyby. {\it The geometry of 
rational surfaces and Hilbert
functions of points in the plane},
Can.\ Math.\ Soc.\ Conf.\ Proc.\ 6 
(1986), 95--111.

\bibitem{\ravello} \manyby. {\it 
Points in Good Position in \pr 2}, in:
Zero-dimensional schemes, Proceedings of the
International Conference held in Ravello, Italy, June 8--13, 1992,
De Gruyter, 1994.

\bibitem{\mtnwest} \manyby. 
{\it Rational surfaces with $K^2>0$}, Proc. A.M.S.
124 (1996), 727--733.

\bibitem{\antican} \manyby. {\it Anticanonical 
rational surfaces}, Trans. A.M.S. 349 (1997), 1191--1208.

\bibitem{\fatpts} \manyby. {\it Free Resolutions of Fat Point 
Ideals on \pr2}, J. Pure and Applied Alg. 125 (1998), 213--234.

\bibitem{\newfatpts} \manyby. {\it The Ideal Generation 
Problem for Fat Points}, preprint, to appear, J. Pure and Applied Alg.

\bibitem{\refHa} Hartshorne, R. Algebraic Geometry. Springer-Verlag, 1977.

\bibitem{\refHi} Hirschowitz, A.
{\it Une conjecture pour la cohomologie 
des diviseurs sur les surfaces rationelles g\'en\'eriques},
Journ.\ Reine Angew.\ Math. 397
(1989), 208--213.

\bibitem{\refMa} Manin, Y.\ I. Cubic Forms. North-Holland 
Mathematical Library
4, 1986.

\bibitem{\refMu} Mumford, D. {\it Varieties defined by quadratic equations},
in: Questions on algebraic varieties, Corso C.I.M.E. 1969 Rome: Cremonese,
1969, 30--100.

\bibitem{\refNtwo} Nagata, M. {\it On rational surfaces, II}, 
Mem.\ Coll.\ Sci.\ 
Univ.\ Kyoto, Ser.\ A Math.\ 33 (1960), 271--293.

\bye